\documentclass[11pt,a4paper]{article}

\usepackage{amsmath, amsfonts, amssymb}
\usepackage{theorem}

\usepackage[latin1]{inputenc}

\newtheorem{thm}{Theorem}[section]

\newtheorem{rem}[thm]{Remark}

\def\be#1\ee{\begin{equation}#1\end{equation}}
\newcommand{\bea}{\begin{eqnarray}}
\newcommand{\eea}{\end{eqnarray}}
\newcommand{\beaa}{\begin{eqnarray*}}
\newcommand{\eeaa}{\end{eqnarray*}}
\newcommand{\bei}{\begin{itemize}}
\newcommand{\eei}{\end{itemize}}
\newcommand{\bee}{\begin{enumerate}}
\newcommand{\eee}{\end{enumerate}}

\def\P{{\mathbb{P}}}

\def\R{\mathbb{R}}

\def\Z{{\mathbb Z}}
\def\N{{\mathbb N}}

\def\ed#1{ {\mathbf 1}_{ \{#1  \}}}             

\def\A{{\bf A}}

\def\aa{{\mathfrak{a}}}

\def\D{{\bf D}}
\def\DD{{\mathfrak D}}
\def\dd{{\rm d}}

\newcommand{\eps}{\varepsilon}
\def\F{{\bf F}}
\def\FF{{\mathfrak F}}

\def\sgn{{\textrm{sgn}}}

\def\i{{\bf i}}

\begin{document}

\title{\bf Small Deviations in $L_2$-norm for Gaussian
Dependent Sequences}
\author{
   Seok Young Hong
   \footnote{Statistical Laboratory, Faculty of Mathematics,
   University of Cambridge,
   UK. email \ {\tt syh30@cam.ac.uk}.
   }
     \and
   Mikhail Lifshits
   \footnote{St.Petersburg State University, Russia,
   and Link\"oping University, Sweden.
   email {\tt mikhail@lifshits.org}.
   }
    \and
    Alexander Nazarov
    \footnote{St.Petersburg Department of Steklov Institute
    of Mathematics
    and St.Petersburg State University, Russia.
    email {\tt al.il.nazarov@gmail.com}.
    }
}
\date{\today}

\maketitle

\begin{abstract}
Let $U=(U_k)_{k\in\Z}$  be a centered Gaussian stationary
sequence satisfying some minor regularity condition. We
study the asymptotic behavior of its weighted $\ell_2$-norm
small deviation probabilities. It is shown that
\[
   \ln \P\left( \sum_{k\in\Z} d_k^2 U_k^2 \le \eps^2\right)
   \sim  - M \eps^{-\frac{2}{2p-1}}, \qquad \textrm{ as }
   \eps\to 0,
\]
whenever
\[
      d_k\sim d_{\pm} |k|^{-p}\quad  \textrm{for some }
      p>\frac{1}{2} \, , \quad k\to \pm\infty,
\]
using the arguments based on the spectral theory of
pseudo-differential operators by M. Birman and M. Solomyak.
The constant $M$ reflects the dependence structure of $U$ in
a non-trivial way, and marks the difference with the
well-studied case of the i.i.d. sequences.
\end{abstract}

\section{Introduction}

Let $(Y(t))_{t\in T}$ be a centered Gaussian process defined
on some parametric measure space $(T,\mu)$. Many studies have
been devoted to the asymptotic behavior of its small deviation
probabilities
\[
   \P\left( ||Y||_2^2= \int\limits_T |Y(t)|^2 \mu(dt)
   \le \eps^2\right), \qquad \textrm{ as } \eps\to 0,
\]
see e.g. \cite{DLL,GH1,GH2,KNN,Naz03,Naz09,NazNik}, to mention
just a small sample. Since by the Karhunen--Lo\`eve expansion
(see for instance \cite[Section 1.4]{AsG})
\[
   ||Y||_2^2 = \sum_{k=1}^\infty d_k^2 X_k^2
\]
where $(X_k)_{k\ge 0}$ is a standard Gaussian i.i.d. sequence
and $d_k^2$ are the eigenvalues of the covariance operator of $Y$,
the small deviation probability may be written as
\[
   \P\left(  \sum_{k=1}^\infty d_k^2 X_k^2 \le \eps^2\right),
   \qquad \textrm{ as } \eps\to 0.
\]
Sharp evaluation of this asymptotics is available when the limiting
behavior of the eigenvalues $d_k^2$ is understood well enough.
Moreover, a considerable amount of results is known also for the
case where $(X_k)$ is an i.i.d.~non-Gaussian sequence, see e.g.
\cite{DLL,Roz4,Roz10}. The importance of small deviation probabilities
in a broader context and the wide spectrum of their applications
are described in the surveys \cite{LiSha,Lif}; for an extensive
up-to-date bibliography see \cite{Bib}.

In this paper, we move towards a different direction and examine
the asymptotic behavior of the small deviation probabilities of
\emph{dependent sequences}. That is,
\be \label{main}
   \P\left(  \sum_{k=1}^\infty d_k^2 U_k^2 \le \eps^2\right),
   \qquad \textrm{ as } \eps\to 0,
\ee
for some stationary centered Gaussian random sequence $U=(U_k)_{k\in\Z}$
that is dependent and only satisfies some mild regularity condition.

The motivation for looking at this small deviation problem under
dependence \eqref{main} is twofold. First, it is an interesting
mathematical question in its own right. The existing literature on
small deviation probability for sums of random variables has been
strictly confined to the i.i.d. framework, so the dependent case is
still an open field of research. Second, there are several potential
statistical applications where this extension could be found useful.
In functional statistics literature, it is well-known that the
convergence rates of nonparametric estimators depend upon the
asymptotics of the associated small deviation probabilities, see e.g.
\cite{FeV}, \cite{Mas} and references therein. Yet in many practical
situations where the functional variable of interest is
discrete-valued, strict independence assumption between the
coordinate variables is too restrictive, so the extent to which the
existing small deviation results can be feasible is limited and the
asymptotics of \eqref{main} should be understood. We refer the
reader to \cite{HoL} for more details.
\medskip

Consider a random vector $Z\in\ell_2(\Z)$ defined by its coordinates
$Z_k= d_k U_k$, $k\in \Z$, where the positive coefficients $d_k$
satisfy the assumption
\be \label{dk}
    d_k\sim d_{\pm} |k|^{-p},\quad
    \textrm{for some } p>\frac{1}{2} \, ,  \quad k\to \pm\infty,
\ee
where at least one of the numbers $d_{\pm}$ is strictly positive.
This assumption is typical of the literature on small deviations of
Gaussian processes and related matters; see for example
\cite{KNN,Li,PW}.

We are interested in the asymptotics of the small deviation
probabilities
\be \label{sd}
     \P\left(||Z||_2 \le \eps\right)
     = \P\left( \sum_{k\in\Z} d_k^2 U_k^2 \le \eps^2\right),
     \qquad \textrm{ as } \eps\to 0.
\ee
In particular, one wonders to what extent this asymptotics is the
same as that for the i.i.d.~Gaussian sequence having the same
variance with $U_k$.

One example of mild dependence structure one can think of would be
{\it linear regularity} (in the sense of \cite[Chapter VII, p.248]{BuS}
and \cite[Chapter 17, p.303]{IbL}). We say that a stationary sequence
$U=(U_k)_{k\in\Z}$ is linearly regular if
\[
   H_{-\infty} := \bigcap_{m\in \Z} H_m = \{0\},
\]
where $H_m$ denotes the closed linear span of $\{ U_k \}_{k\le m}$.
It is a type of asymptotic independence condition that roughly means
the process has no significant influence from the distant past. When
the process is Gaussian, linear regularity is implied by the class of
mixing-type conditions, a popular notion of dependence under which
probability theories have been extensively studied in the literature;
see e.g. \cite{Br} and \cite{Doukh} for the precise definition and a
comprehensive review.

Since a consequence of the Wold decomposition theorem suggests that
any stationary linearly regular Gaussian sequence admits a causal
moving average representation
(cf. \cite[Chapter VII, Theorem 13]{BuS}):
\[
   U_k = \sum_{m=0}^\infty a_m X_{k-m}
   =  \sum_{j=-\infty}^k a_{k-j} X_{j} ,
\]
where $\sum_{m=0}^\infty a_m^2<\infty$ and $(X_j)_{j\in\Z}$ is an
i.i.d.~standard Gaussian sequence, it follows that many popular
dependent processes  such as strongly mixing sequences
do have such representations.

In the sequel we shall consider a more general assumption than
causality, and postulate that
\be \label{Uk}
   U_k = \sum_{m=-\infty}^\infty a_m X_{k-m}
   =  \sum_{j=-\infty}^\infty a_{k-j} X_{j} ,
\ee
where $(a_m)\in \ell_2(\Z)$, and $(X_j)$ is i.i.d. standard Gaussian
as above. In fact, this representation exists iff the stationary sequence
$(U_k)$ has a spectral density (cf. Remark \ref{spd} below) but we will
not develop this point of view any further.

Our main result is as follows:

\begin{thm} \label{t:sd_statseq}
Let a stationary centered Gaussian sequence $(U_k)_{k\in\Z}$ admit a
representation \eqref{Uk} and let the coefficients $(d_k)_{k\in\Z}$
have the asymptotics \eqref{dk}. For $p<1$ suppose in addition that
$(a_m)\in \ell_r(\Z)$ with some $r<2$. Then
\be \label{sd2}
   \ln \P\left( \sum_{k\in\Z} d_k^2 U_k^2 \le \eps^2\right)
   \sim - B_p \,\left(\frac{C}{\eps^2}\right)^{\frac 1{2p-1}},
   \qquad \textrm{ as } \eps\to 0,
\ee
with the constants
\[
   B_p=  \frac{2p-1}{2}\left(\frac{\pi}{2p\sin
   \big(\frac{\pi}{2p}\big)}\right)^{\frac{2p}{2p-1}},
\]
\be\label{C}
   C= \left( \frac{1}{2\pi} \int\limits_0^{2\pi}
      \Bigl|\sum_{m=-\infty}^\infty
      a_m\, e^{\i\, m x}\Bigr|^{1/p} dx\right)^{2p}
     \left( d_-^{1/p}+ d_+^{1/p}\right)^{2p}.
\ee
\end{thm}

 \begin{rem}
 {\rm
The power term in the logarithmic small deviation asymptotics is the
same as that in the i.i.d. case (characterized by $a_m=a_0 \ed{m=0}$),
but the constant $C$ in front of it depends on the sequence
$(a_m)
$ in a nontrivial way, no matter how weak the linear dependence
in $(U_k)$ is (in other words, how fast $a_m$ decays).
}
\end{rem}

\begin{rem}
{\rm We do not know whether the extra assumption on $(a_m)$ for
$p<1$ is essential or purely technical.}
\end{rem}

\begin{rem}
{\rm For sharper results on small deviations, one would need
to know a sharper spectral asymptotics (the so-called two-term
asymptotics). This seems to be a much harder problem in general.}
\end{rem}

\begin{rem}
{\rm Similar technique can be applied in the study of the weighted
$L_2$-norm small deviations for {\it continuous time} stationary
processes. This will be done elsewhere.}
\end{rem}

\section{Proof of Theorem \ref{t:sd_statseq}}

Recall that we have a random vector $Z=(d_k U_k)\in \ell_2(\Z)$ and
a random vector with independent coordinates $X=(X_j)$, $j\in\Z$.
It follows from the definitions that
\[
    Z= \D U = \D \A X,
\]
where $\D$ is the diagonal matrix with elements $d_{kj}=d_k \ed{k=j}$
and $\A$ is the Toeplitz matrix with elements $a_{kj}= a_{k-j}$.
Therefore, the covariance operator of $Z$ that maps $\ell_2(\Z)$ into
$\ell_2(\Z)$ can be expressed as
\[
    K_Z = cov(Z)= (\D \A) (\A^* \D),
\]
and by the Karhunen--Lo\`eve expansion (see \cite[Section 1.4]{AsG}),
\[
    ||Z||^2= \sum_{n=1}^\infty \lambda_n \, \xi_n^2\ ,
\]
where $(\xi_n)_{n\in\N}$ is an i.i.d. standard Gaussian sequence and
$(\lambda_n)_{n\in\N}$ are the eigenvalues of $K_Z$.

We remark that the small deviations \eqref{sd} depend heavily on the
asymptotic behavior of $\lambda_n$. In particular, if we can show that
\be \label{lambda}
     \lambda_n  \sim  C\, n^{-2p},
     \qquad \textrm{ as } n\to\infty,
\ee
then \eqref{sd2} will follow from  \cite[p.67]{DLL} or \cite{Zol}, and
\cite{Naz09}. The decay rate for $\lambda_n$ would then be the same as
that of $d_n^2$, and the constant $C$ in front of the power rate would
depend on the sequence $(a_m)$ in a non-trivial way, cf. \eqref{C}.

Therefore it now remains to prove the eigenvalue asymptotics \eqref{lambda},
and to specify the constant $C$.\\
%

Since all separable Hilbert spaces are isomorphic, we may replace
$\ell_2(\Z)$ with the more appropriate space
$L_2\left([0,2\pi],\nu \right)$ with $\nu(dy)=\tfrac{dy}{2\pi}$,\,
equipped with the standard exponential basis
$e_m(x)=\exp( \i\, m x),~m\in \Z$.

Notice that in this space $\A$ becomes the multiplication operator
$\A f = \aa f$ related to the function
\[
   \aa(x)=\sum_{m=-\infty}^\infty  a_m\, e_m(x),
\]
while $\D$ becomes the convolution operator
\[
   (\D f)(x) =\int\limits_0^{2\pi} \DD(x-y) f(y)\, \nu(dy)
\]
with the kernel
\be \label{D}
    \DD(x) =\sum_k d_k\, e_k(x).
\ee
Indeed, if $f=\sum_j f_j \, e_j$, then
\[
   \aa f = \sum_{m,j} a_m f_j \, e_{m+j}
   = \sum_k \left(\sum_j a_{k-j} f_j \right) e_k
\]
and
\begin{eqnarray*}
  \int\limits_0^{2\pi}  \DD(x-y) f(y)  \, \nu(dy)
  &=&  \sum_{j,k}  d_k f_j
       \int\limits_0^{2\pi}  e_k(x-y) e_j(y) \,\nu(dy)
  \\
  &=& \sum_{j,k} d_k f_j  e_k(x)
      \int\limits_0^{2\pi}  e_{j-k}(y) \,\nu(dy)
  =  \sum_{k} d_k f_k  e_k(x).
\end{eqnarray*}

\begin{rem} \label{spd}
Interestingly, $|\aa(\cdot)|^2$ is the spectral density
of the stationary sequence $(U_k)$.
\end{rem}
\medskip

In our spectral analysis, we will first slightly reinforce
condition \eqref{dk} by assuming that $(d_k)$ is {\it exactly}
equal to the non-isotropic power function
\be \label{homogen}
  d_k = d(\sgn(k)) \  |k|^{-p},
\ee
where $d(\pm 1)=d_{\pm}$ are two constants and $d_0=0$.

In the sequel, our main argument will be a reduction of the
operator $\A^* \D$ to a special case of the pseudo-differential
operators ($\Psi$DO) studied by M. Birman and M. Solomyak
(hereafter BS) in \cite{BS1,BS2}\footnote{The referee
mentioned \cite{BS3} which also provides estimates relevant
to small deviations. However, these estimates are not sharp
enough to establish the asymptotic behavior of singular
values up to equivalence that we need here.},
see also \cite{DauRob}.

The following exposition provides an interpretation of \cite{BS1}
and \cite{BS2} adapted to our case. The aim of the papers BS is
the spectral analysis of the following operator (in their notation)
\[
   (\F u)(x) = b(x)\int\limits_{\R^m} \FF(x,x-y) c(y) u(y) dy.
\]
Here and elsewhere by spectral analysis of an operator, we
understand the study of the asymptotic behavior of its singular
values.

In our case the space dimension $m=1$, and we can assume that the
function $\FF$ depends only on the second argument, i.e.
\be \label{bs1a}
  (\F u)(x) = b(x) \int\limits_{\R} \FF(x-y) c(y) u(y) dy.
\ee

The kernel $\FF(\cdot)$ in \cite{BS1} is of specific Fourier
transform form, namely,
\be \label{bs_F}
  \FF =
  (\zeta\cdot \dd )\check{}.
\ee
Here $\zeta(\cdot)$ is any smooth function that vanishes on
a neighborhood of zero and equals to one on a neighborhood of
infinity, while $\dd(\cdot)$ in the one-dimensional case is a
homogeneous function as in \eqref{homogen} but considered in
continuous time, i.e., in the notation of BS
\be \label{homogen2}
  \dd(\xi) = \dd(\sgn(\xi)) \  |\xi|^{-\alpha},
  \qquad \xi\in\ \R\backslash\{0\},
\ee
where $\dd(\pm 1)=\dd_{\pm}$ are two constants. For us,
the homogeneity index $\alpha$ in \eqref{homogen2} is $p$.
Notice immediately that the ``mysterious'' formula
\eqref{bs_F} is, apart from the inessential smoothing term
$\zeta$, a version of our former kernel definition \eqref{D}
for continuous time.

BS consider the operator $\F$ either on $\R^m$ or on a cube.
The latter means that the weights $b$ and $c$ in \eqref{bs1a}
are supported by a cube. In our case the weight function
$b(\cdot)$ from \eqref{bs1a} is $\overline{\aa(\cdot)}$, and
the function $c(\cdot)$ is the indicator on the interval
$[0,2\pi]$ that plays the role of a cube. Moreover, the index
$\mu=\tfrac{m}{\alpha}$ used by BS for the description of
singular values behavior is $\tfrac{1}{p}$ in our notation.
Notice that \cite{BS1} distinguishes three cases $\mu>1$,
$\mu=1$ and $\mu<1$, which in our notation are
$p\in(\frac 12, 1)$, $p=1$ and  $p>1$, respectively.

The weight size restrictions in \cite{BS1} are $b\in L_{q_1}$,
$c\in L_{q_2}$. Our assumptions give $q_1=2$ for $p\ge 1$ and
$q_1=\tfrac{r}{r-1}>2$ for $p<1$ (the latter fact is due to
the Hausdorff--Young inequality, see, e.g. \cite[\S\,8.5]{HLP}).
Without loss of generality we can suppose
$\tfrac{r}{r-1}<\frac 1p$. Further, $q_2\ge 1$ may be taken
arbitrarily.
\medskip

The main results of BS are stated in Theorems 1 and 2 of
\cite{BS1}. Let us first check the weight assumptions of
Theorem 1 in \cite{BS1}.

If $p>1$, then $\mu<1$ and Theorem 1(b) applies with
$q_1=q_2=2$.

If $p=1$, then $\mu=1$ and Theorem 1(c) applies with $q_1=2$
and any $q_2>2$. This case is relevant to Wiener process and
its relatives such as Brownian bridge, OU-process etc.

If $p\in(\tfrac 12,1)$, then $2>\mu>1$, and Theorem 1(a)
applies with $q_1>2$ and $q_2>2$ chosen from the relation
$\tfrac{1}{q_2} = p-\tfrac{r-1}{r}$, as required in
Theorem 1(a).

Theorem 2 in \cite{BS1} is disregarded because it
requires some extra assumptions and only applies to the case of
infinite $q_1$ or $q_2$.
\medskip

Now let us proceed to follow the BS result. They denote the
singular values of $\F$ by $s_n(\F)$ and study the corresponding
{\it distribution function}
\[
  N_{\F}(s) :=  \#\{n: s_n(\F)\ge s\}
\]
and its asymptotics at zero. This is indeed an equivalent setting
because
\be \label{Nass}
    N_{\F}(s)  \sim \Delta\cdot s^{-1/p},
    \quad \textrm{as } s\to 0 \quad
    \Longleftrightarrow \quad
    s_n(\F) \sim \Delta^p \cdot n^{-p},
    \quad \textrm{as } n\to \infty.
\ee

Next, BS introduce
the following notations
\be
   \Delta_\mu := \limsup_{s\to 0_+} s^\mu N_{\F}(s), \qquad
   \delta_\mu := \liminf_{s\to 0_+} s^\mu N_{\F}(s).
\ee

In their Theorem 2 of \cite{BS1} BS prove that
$\Delta_\mu= \delta_\mu$ and
find
the common value for the upper and the lower limit
\[
   \lim_{s\to 0_+} s^\mu N_{\F}(s)= \Delta_\mu= \delta_\mu.
\]
Namely, they introduce the ``operator symbol''
$G(s,\xi)$, see formula (14) of \cite{BS1}. In the one-dimensional
case the symbol is a scalar defined by
\[
   G(x,\xi) = \overline\aa(x) \dd(\xi) =\overline\aa(x)\cdot
   {\mathbf 1}_{[0,2\pi]}(x)\cdot \dd(\sgn(\xi)) \, |\xi|^{-p}.
\]

Further, formula (18) of \cite{BS1} suggests that in our case (recall
that $\mu=\tfrac{1}{p}$)
\begin{eqnarray*}
   \Delta_\mu &=&  (2\pi)^{-1}  \int\limits_0^{2\pi}
   \int\limits_{\R\backslash\{0\}} \ed{|G(x,\xi)|\ge 1}  d\xi dx
   \\
   &=&  (2\pi)^{-1} \int\limits_0^{2\pi} \int\limits_{\R\backslash\{0\}}
       \ed{|\aa(x)|\, |\dd(\sgn(\xi))| \, |\xi|^{-p} \ge 1}  d\xi dx
   \\
   &=&  (2\pi)^{-1} \int\limits_0^{2\pi} \int\limits_{\R\backslash\{0\}}
        \ed{ |\aa(x)|^{1/p}\, |\dd(\sgn(\xi))|^{1/p} \, \ge |\xi|}  d\xi dx
   \\
   &=&  (2\pi)^{-1} \int\limits_0^{2\pi}  |\aa(x)|^{1/p} dx
       \left( |\dd(-1)|^{1/p}+ |\dd(1)|^{1/p} \right).
\end{eqnarray*}

Now we compare the spectral behavior of the operator of our interest
$\A^* \D$ with that of the operator $\F$ in \eqref{bs1a}, assuming
that the parameters $d_{\pm}$ in \eqref{homogen} coincide with their
counterparts
$\dd_{\pm}$ in \eqref{homogen2}, and substituting
$b=\overline\aa\cdot {\mathbf 1}_{[0,2\pi]}$ and
$c= {\mathbf 1}_{[0,2\pi]}$ in \eqref{bs1a}.

Let us prove that
\be\label{equal}
      N_{\A^* \D}(s)\sim N_{\F}(s), \quad \textrm{as } s\to 0.
\ee

Notice that since we are working on the interval of length $2\pi$,
it is sufficient to consider only the restriction of our periodical
function $\DD$ to $[-2\pi,2\pi]$.

Let $h$ be the cut-off function equal to one on $[\tfrac{3\pi}{2},2\pi]$
and zero on $[-2\pi,\pi]$. Then it follows that the function
$h_0(x):=1- h(x)-h(-x)$ equals to one on $[-\pi,\pi]$ and vanishes outside
of the interval  $[-\tfrac{3\pi}{2},\tfrac{3\pi}{2}]$.

Comparing the kernels of two operators, we have the following decomposition
\be  \label{decomp}
   \DD(x)-\FF(x) = \DD(x)\ \big(h(x)+h(-x)\big) + \DD_1(x),
   \qquad x\in [-2\pi,2\pi].
\ee

We claim that the function $\DD_1:= \DD\cdot h_0(x)-\FF$ satisfies
\be \label{op}
   \widehat{\DD_1}(\xi)=o(|\xi|^{-p})\quad \textrm{as } |\xi|\to\infty,
\ee
where $\widehat{\DD_1}$ denotes the Fourier transform of $\DD_1$.
Indeed, we have
\[
   \widehat{\DD\cdot h_0}(\xi)=\sum_{k\ne0} d(\sgn(k)) \
   |k|^{-p}\widehat h_0(\xi-k),
\]
and then by spliting the series into two sums,
\[
  \widehat{\DD\cdot h_0}(\xi)=\Sigma_1+\Sigma_2
  :=\left(\sum_{|k-\xi|\le\sqrt{\xi}}+\sum_{|k-\xi|>\sqrt{\xi}}\right)
  d(\sgn(k)) \  |k|^{-p}\widehat h_0(\xi-k).
\]
Since $\widehat h_0$ rapidly decays at infinity, we have
$\Sigma_2=o(|\xi|^{-p})$ as $|\xi|\to\infty$. Further,
\begin{eqnarray*}
  \Sigma_1 &=& d(\sgn(\xi)) \  |\xi|^{-p}
  \sum_{|k-\xi|\le\sqrt{\xi}}\widehat h_0(\xi-k)+o(|\xi|^{-p})
  \\
  &=& d(\sgn(\xi)) \  |\xi|^{-p}\sum_k\widehat h_0(\xi-k)
     +o(|\xi|^{-p})=d(\sgn(\xi)) \  |\xi|^{-p}+o(|\xi|^{-p})
\end{eqnarray*}
by the Poisson summation formula (see, e.g.,
\cite[Ch. II, Sect. 13]{Zyg}), so that (\ref{op}) follows.

Decomposition \eqref{decomp} generates the corresponding operator
representation
\[
   \A^*\D -\F = (\D_+  +  \D_-) + \D_1.
\]

By corollary 4) in \cite{BS2}, relation \eqref{op} gives
$\lim_{s\to 0_+} s^{1/p} N_{\D_1}(s)=0$. Further, since
$\DD$ is $2\pi$-periodic, the singular values of $\D_+$ coincide
with the singular values of the operator
\[
    \overline{\aa}(x+\pi) {\mathbf 1}_{[0,\pi]}(x)
    \int\limits_{\R}
    \DD(x-y)h(x+2\pi-y) {\mathbf 1}_{[\pi,2\pi]}(y) u(y) dy.
 \]
For this operator, we have $\rm{supp}(b)=[0,\pi]$ and
$\rm{supp}(c)=[\pi, 2\pi]$ in terms of (\ref{bs1a}),
and Lemma 3 in \cite{BS2} gives
$\lim_{s\to 0_+} s^{1/p} N_{\D_+}(s)=0$. By the same reason,
$\lim_{s\to 0_+} s^{1/p} N_{\D_-}(s)=0$, yielding \eqref{equal}.
\medskip

Using the equivalence in \eqref{Nass}, we obtain
\begin{eqnarray*}
   s_n(\A^* \D) &\sim& \Delta_\mu^p \, n^{-p}
\\
  &=& \left( \frac{1}{2\pi}
      \int\limits_0^{2\pi} |\aa(x)|^{1/p} dx\right)^p
      \left( |d(-1)|^{1/p} + |d(1)|^{1/p}\right)^p n^{-p}.
\end{eqnarray*}
Since $\lambda_n=s_n^2(\A^*\D)$ by the definition of singular values,
it follows that
\[
     \lambda_n  \sim
        \left( \frac{1}{2\pi} \int\limits_0^{2\pi}
        |\aa(x)|^{1/p} dx\right)^{2p}
        \left( |d(-1)|^{1/p}+ |d(1)|^{1/p}\right)^{2p} n^{-2p},
        \qquad n\to\infty,
\]
as required in \eqref{lambda}, and the conclusion for small
deviations follows.
\medskip

So far, the result of the theorem is obtained  only for the
{\it homogeneous} coefficients \eqref{homogen}. However, since
any finite number of terms in the sequence $(d_k)$ is irrelevant
for small deviation probability asymptotics, by monotonicity of
the quadratic form $\sum_{k\in\Z} d_k^2 U_k^2$ in $(d_k)$, it
follows that \eqref{sd2} also holds for
any $(d_k)$ satisfying \eqref{dk}. \hfill $\Box$
\medskip

{\bf Acknowledgment.} We are grateful to Professor G. Rozenblum
for useful advice and for a number of provided references and
to the anonymous referee for careful reading, constructive
critics and for a number of suggested improvements.

M.A.\,Lifshits and A.I.\,Nazarov are supported by RFBR grant
16-01-00258. S.Y.\,Hong acknowledges financial support from the
ERC grant 2008-AdG-NAMSEF.


\begin{thebibliography}{1}

{\baselineskip=12pt \parskip=2pt

\bibitem{AsG}
   R.B.\,Ash, M.F.\,Gardner. \emph{Topics in Stochastic Processes}.
   Academic Press, New York, 1975.

\bibitem{BS1}
   M.\v{S}.\,Birman, M.Z.\,Solomjak. {\it Asymptotics of the spectrum
   of pseudodifferential operators with anisotropic-homogeneous
   symbols}, Vestnik LGU (1977), no 13, 13--21 (Russian); English
   transl.: Vestnik Leningrad Univ. Math. \textbf{10} (1982),
   237--247.

\bibitem{BS2}
   M.\v{S}.\,Birman, M.Z.\,Solomjak. \emph{Asymptotics of the
   spectrum of pseudodifferential operators with
   anisotropic-homogeneous symbols. II}, Vestnik LGU (1979),
   no~13, 5--10 (Russian); English transl.: Vestnik Leningrad Univ.
   Math. \textbf{12} (1980), 155--161.

\bibitem{BS3}
   M.\v{S}.\,Birman, M.Z.\,Solomjak. \emph{Estimates of singular
   numbers of integral operators}, Uspekhi Mat. Nauk \textbf{32}
   (1977), no 1, 17--84 (Russian); English transl.: Russian Math.
   Surveys \textbf{32} (1977), 15--89.

\bibitem{Br}
   R.C.\,Bradley. \emph{Introduction to Strong Mixing Conditions},
   Vol. 1--3. Kendrick Press, Heber City, 2007.

\bibitem{BuS}
   A.V.\,Bulinskii, A.N.\,Shiryaev. \emph{Theory of Random
   Processes}, Fizmatlit, 2003 (in Russian).

\bibitem{DauRob}
   M.\,Dauge, D.\,Robert. \emph{Weyl's formula for a class of
   pseudodifferential operators with negative order on
   $L_2(\R^n)$.} In: {Proc. Conf. ``Pseudo-differential
   operators'', Oberwolfach, 1986}; Lecture Notes in Math.,
   v.~1256, Springer-Verlag, 1987, 90--122.

\bibitem{Doukh}
   P.\,Doukhan. \emph{Mixing: Properties and Examples.} Lecture
   Notes in Statistics 85, Springer-Verlag, 1994.

\bibitem{DLL}
   T.\,Dunker,  M.A.\,Lifshits,  W.\,Linde. \emph{Small deviations
   of sums of independent variables.} In: {Proc. Conf. High
   Dimensional Probability}; Ser. Progress in Probability,
   v.~43, Birkh\"auser, 1998, 59--74.

\bibitem{FeV}
   F.\,Ferraty, P.\,Vieu. \emph{Nonparametric Functional Data
   Analysis: Theory and Practice}, Springer, New York, 2006.

\bibitem{GH1}
   F.\,Gao, J.\,Hannig, T.-Y.\,Lee, and F.\,Torcaso. \emph{Laplace
   transforms via {H}adamard factorization with applications to
   small ball probabilities}, Electronic J. Probab. \textbf{8}
   (2003), paper 13.

\bibitem{GH2}
   F.\,Gao, J.\,Hannig, T.-Y.\,Lee, and F.\,Torcaso. \emph{Exact
   {$L^2$}-small balls of Gaussian processes}, J. Theoret. Probab.
   \textbf{17} (2004), no.~2, 503--520.

\bibitem{HLP}
   G.H.\,Hardy, J.E.\,Littlewood, G. P\'olya. \emph{Inequalities},
   Cambridge University Press, Cambridge, 1934.

\bibitem{HoL}
  S.Y.\,Hong, O.\,Linton. \emph{Asymptotic properties of a
  Nadaraya--Watson type estimator for regression functions of
  infinite order}.  Preprint {\tt https://arxiv.org/abs/1604.06380}.

\bibitem{IbL}
  I.A.\,Ibragimov, Yu.V.\,Linnik. \emph{Independent and Stationary
  Sequences of Random Variables}, Wolters-Noordhoff, Groningen,
  1971.

\bibitem{KNN}
  A.I.\,Karol', A.I.\,Nazarov, Ya.Yu.\,Nikitin. \emph{Small ball
  probabilities for {G}aussian random fields and tensor products
  of compact operators}, Trans. Amer. Math. Soc. \textbf{360}
  (2008), no.~3, 1443--1474.

\bibitem{Li}
  W.V.\,Li. \emph{Comparison results for the lower tail of
  {G}aussian seminorms}, J. Theor. Probab. \textbf{5} (1992),
  1--31.

\bibitem{LiSha}
  W.V.\,Li, Q.-M.\,Shao. \emph{Gaussian processes: inequalities,
  small ball probabilities and applications}, In: Stochastic
  Processes: Theory and Methods, Handbook of Statistics (C.R. Rao
  and D.~Shanbhag, eds.), vol.~19, North-Holland/Elsevier,
  Amsterdam, 2001, pp.~533--597.

\bibitem{Lif}
  M.A.\,Lifshits. \emph{Asymptotic behavior of small ball
  probabilities}, In: Probab. Theory and Math. Statist. Proc.
  VII International Vilnius Conference (1998) (B.~Grigelionis, ed.),
  VSP/TEV. Vilnius, 1999, pp.~453--468.

\bibitem{Bib}
  M.A.\,Lifshits. \emph{Bibliography of small deviation
  probabilities},\quad On the small deviation website
  {\tt http://www.proba.jussieu.fr/pageperso/smalldev/biblio.pdf}

\bibitem{Mas}
  A.\,Mas. \emph{Lower bound in regression for functional data by
  small ball probability representation in Hilbert space},
  Electronic J. Statist. \textbf{6} (2012), 1745--1778.

\bibitem{Naz03}
  A.I.\,Nazarov. \emph{Exact {$L_2$}-small ball asymptotics of
  Gaussian processes and the spectrum of boundary-value problems},
  J.\,Theor.\,Probab. \textbf{22} (2009), no.~3, 640--665.

\bibitem{Naz09}
  A.I.\,Nazarov. \emph{Log-level comparison principle for small ball
  probabilities}. Statist.\,\&\,Probab.\,Letters \textbf{79} (2009),
  no.~4, 481--486.

\bibitem{NazNik}
  A.I.\,Nazarov, Ya.Yu.\,Nikitin. \emph{Exact $L_2$-small ball behavior
  of integrated Gaussian processes and spectral asymptotics of boundary
  value problems}. Probab.\,Theor.\,Rel.\,Fields \textbf{129} (2004), 
  no.~4, 469--494.

\bibitem{PW}
  A.\,Papageorgiou, G.W.\,Wasilkowski. \emph{On the average complexity
  of multivariate problems},  J.\,Complexity \textbf{6} (1990), 1--23.

\bibitem{Roz4}
L.V.\,Rozovsky. \emph{Small deviation probabilities for sums of independent
  positive random variables}, J.\,Math.\,Sci. \textbf{147} (2007), no.~4,
  6935--6945.

\bibitem{Roz10}
  L.V.\,Rozovsky. \emph{On the behavior of the log Laplace transform of
  series of weighted non-negative random variables at infinity},
  Statist. \& Probab. Letters \textbf{80} (2010), 764--770.

\bibitem{Zol}
  V.M.\,Zolotarev. \emph{Asymptotic behavior of Gaussian measure in
  $\ell_2$}, J. Math. Sci. \textbf{35} (1986), 2330--2334.

\bibitem{Zyg}
  A.\,Zygmund. \emph{Trigonometrical Series}, Vol.1, Cambridge University
  Press, Cambridge, 1959.

} 

\end{thebibliography}
\end{document}